\documentclass[12pt, a4paper]{article}
\usepackage[latin1]{inputenc}
\usepackage[english]{babel}
\usepackage{amsmath}
\usepackage[english]{babel}
\usepackage{indentfirst}
\usepackage{amstext}
\usepackage{enumerate}
\usepackage{amsfonts}
\usepackage{textcomp}
\usepackage{amssymb}
\usepackage[dvips]{graphicx}
\usepackage{setspace}
\usepackage[all]{xy}
\usepackage{color}

\newcommand {\demo}{\hskip -0.6cm{\bf Proof.  }}
\newcommand {\fim}{\hfill{$\square$}\vskip 1pc}

\newcommand {\N}{\mathbb{N}}

\newcommand {\GG}{\mathcal{G}}
\newcommand {\pinf}{\mathfrak{p}^\infty}
\newcommand {\pfin}{\mathfrak{p}^*}
\newtheorem{teorema}{Theorem}[section]
\newtheorem{lema}[teorema]{Lemma}
\newtheorem{corolario}[teorema]{Corollary}
\newtheorem{definicao}[teorema]{Definition}
\newtheorem{proposicao}[teorema]{Proposition}
\newtheorem{exemplo}[teorema]{Example}
\newtheorem{remark}[teorema]{Remark}

\begin{document}
\onehalfspace

\title{Irreducible and permutative representations of ultragraph Leavitt path algebras}

\author{Daniel Gon\c{c}alves\footnote{This author is partially supported by Conselho Nacional de Desenvolvimento Cient\'{i}fico e Tecnol\'{o}gico - CNPq and Capes-PrInt.}  \ and Danilo Royer}
\maketitle

\begin{abstract}
  We completely characterize perfect, permutative, irreducible representations of an ultragraph Leavitt path algebra. For this we extend to ultragraph Leavitt path algebras Chen's construction of irreducible representations of Leavitt path algebras. We show that these representations can be built from branching system and characterize irreducible representations associated to perfect branching systems. Along the way we improve the characterization of faithfulness of Chen's irreducible representations.
\end{abstract}

\vspace{1.0pc}
MSC 2010: 16G99, 16S10, 16W50.

\vspace{1.0pc}
Keywords: Ultragraph Leavitt path algebras, irreducible representations, permutative representations, branching systems.

\section{Introduction}

Ultragraphs (and their algebras) are generalizations of graphs (and their algebras) with applications in Symbolic Dynamics, Operator Algebras and Algebra, see for example \cite{GD, GLR1, dd3, dd1, GRultrapartial, GRultra, GSCSC, LiYorke, brunodaniel,  KMST, afultra, larki, Tomforde:JOT03, Tomsimplicity}. The Leavitt path algebra associated to an ultragraph was defined in \cite{leavittultragraph}, where it is shown that these algebras provide examples of algebras that can not be realized as the Leavitt path algebra of a graph. A key feature of an ultragraph path algebra is that it provides a unified approach to both Leavitt path algebras and Cuntz-Krieger algebras associated with infinite matrices (see \cite{leavittultragraph} for the purely algebraic context and \cite{Tomforde:JOT03} in the C*-algebraic context). It is therefore interesting to extend known results of Leavitt path algebra theory to ultragraph Leavitt path algebras. Furthermore, since ultragraph Leavitt path algebras are algebraic analogues of ultragraph C*-algebras, which are well studied and play a key role in the study of infinite alphabet shift spaces (see \cite{GRultrapartial, GRultra, GSCSC, LiYorke, brunodaniel}), it is important to deepen the understanding of these algebras.

Recently there has been intense activity on the study of representations of Leavitt path algebras. For example, in \cite{colorado, malaga, malaga1} it is shown that  irreducible representations play an important role in the study of the socle series of Leavitt path algebras. The study of representations via branching systems was done in \cite{GR} and in \cite{chen} a key type of irreducible representation was constructed (which is now called a Chen module). The investigation of Leavitt path algebras with a special type of, or a specific number of, irreducible representations was done in  \cite{A6, colorado1, AR, HR, ranga, Ranga}.

Our goal in this paper is to contribute to the study of representations of ultragraph Leavitt path algebras. In particular we will extend Chen's results (see \cite{chen}) regarding irreducible representations of Leavitt path algebras to ultragraph algebras, improve some of them, and use our results to describe permutative, perfect, irreducible representations (a result that is new also in the context of graph algebras). Our interest in permutative representations come from the fact that they have applications to wavelets, continued fraction expansions, iterated function systems, higher rank graphs, among others, see \cite{bratteli, MR3404559, kawamura8}.

The paper is organized as follows: after this introduction we include a brief sections of preliminaries, which is followed by Section~3, where we extend Chen's representations to ultragraph Leavitt path algebras. In 
Section~4 we show that the representations build in Section~3 can be obtained via branching systems, and use branching system theory to completely characterize faithfulness of the representations (this result improves known results for Leavitt path algebras of graphs). We focus on perfect representations and perfect branching systems in Section~5, where we completely characterize irreducible representations of ultragraph path algebras arising from perfect branching system as those of Section~3 (this extends results of \cite{chen}). Finally, in Section~6, we completely characterize perfect, permutative, irreducible representations of an ultragraph Leavitt path algebra (this is a new result also in the context of Leavitt path algebras of graphs).

\section{Preliminaries}

In this brief section we recall the definition of the Leavitt path algebra associated to an ultragraph and set notation. In particular, unless otherwise stated, we let $R$ denote a unital commutative ring throughout the paper.

\begin{definicao}\label{def of ultragraph}
An \emph{ultragraph} is a quadruple $\mathcal{G}=(G^0, \mathcal{G}^1, r,s)$ consisting of two countable sets $G^0, \mathcal{G}^1$, a map $s:\mathcal{G}^1 \to G^0$, and a map $r:\mathcal{G}^1 \to P(G^0)\setminus \{\emptyset\}$, where $P(G^0)$ stands for the power set of $G^0$.
\end{definicao}

\begin{definicao}\label{def of mathcal{G}^0}
Let $\mathcal{G}$ be an ultragraph. Define $\mathcal{G}^0$ to be the smallest subset of $P(G^0)$ that contains $\{v\}$ for all $v\in G^0$, contains $r(e)$ for all $e\in \mathcal{G}^1$, and is closed under finite unions and non-empty finite intersections. Elements of $\mathcal{G}^0$ are called generalized vertices.
\end{definicao}


\begin{definicao}\label{def of ultragraph algebra}
Let $\mathcal{G}$ be an ultragraph and $R$ be a unital commutative ring. The Leavitt path algebra of $\mathcal{G}$, denoted by $L_R(\mathcal{G})$, is the universal algebra with generators $\{s_e,s_e^*:e\in \mathcal{G}^1\}\cup\{p_A:A\in \mathcal{G}^0\}$ and relations
\begin{enumerate}
\item $p_\emptyset=0,  p_Ap_B=p_{A\cap B},  p_{A\cup B}=p_A+p_B-p_{A\cap B}$, for all $A,B\in \mathcal{G}^0$;
\item $p_{s(e)}s_e=s_ep_{r(e)}=s_e$ and $p_{r(e)}s_e^*=s_e^*p_{s(e)}=s_e^*$ for each $e\in \mathcal{G}^1$;
\item $s_e^*s_f=\delta_{e,f}p_{r(e)}$ for all $e,f\in \mathcal{G}$;
\item $p_v=\sum\limits_{s(e)=v}s_es_e^*$ whenever $0<\vert s^{-1}(v)\vert< \infty$.
\end{enumerate}
\end{definicao}
Let $\mathcal{G}$ be an ultragraph. A finite path is either an element of $\mathcal{G}^0$ or a sequence of edges $e_1...e_n$, with length $|e_1...e_n|=n$, and such that $s(e_{i+1})\in r(e_i)$ for each $i\in \{0,...,n-1\}$. An infinite path is a sequence $e_1e_2e_3...$, with length $|e_1e_2...|=\infty$, such that $s(e_{i+1})\in r(e_i)$ for each $i\geq 0$. The set of finite paths in $\mathcal{G}$ is denoted by $\mathcal{G}^*$, and the set of infinite paths in $\mathcal{G}$ is denoted by $\mathfrak{p}^\infty$. We extend the source and range maps as follows: $r(\alpha)=r(\alpha_{|\alpha|})$, $s(\alpha)=s(\alpha_1)$ for $\alpha\in \mathcal{G}^*$ with $0<|\alpha|<\infty$, $s(\alpha)=s(\alpha_1)$ for each $\alpha\in \mathfrak{p}^\infty$, and $r(A)=A=s(A)$ for each $A\in \mathcal{G}^0$. An element $v\in G^0$ is a sink if $s^{-1}(v) = \emptyset$, and we denote the set of sinks in $G^0$ by $G^0_s$. We say that $A\in \mathcal{G}^0$ is a sink if each vertex in $A$ is a sink. We also define the set $\mathfrak{p}^*$ by 
$$\pfin=\{(\alpha,v): \alpha\in \GG^*, |\alpha|\geq 1 \text{ and }v\in r(\alpha)\cap G_s^0\}\cup \{(v,v): v \in G^0_s\}.$$

\begin{remark}
Notice that given a vertex $v$, the element $(v,v)$ is an element of $\mathfrak{p}^*$ if, and only if, $v$ is a sink. 
\end{remark}

\begin{exemplo}\label{ex1}
Let $\GG$ be the ultragraph as follows:

\vspace{2cm}

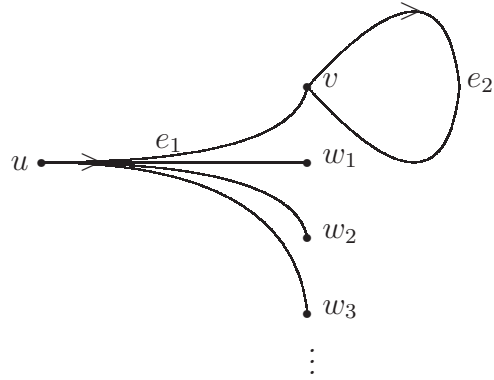
\centerline{
\setlength{\unitlength}{1cm}
\begin{picture}(4,0)
\put(-0.5,0){\circle*{0.1}}
\put(-0.9,-0.1){$u$}
\qbezier(-0.5,0)(2.8,0)(3,1)
\qbezier(-0.5,0)(2.8,0)(3,0)
\qbezier(-0.5,0)(2.8,0)(3,-1)
\qbezier(-0.5,0)(2.8,0)(3,-2)
\put(3,1){\circle*{0.1}}
\put(3,0){\circle*{0.1}}
\put(3,-1){\circle*{0.1}}
\put(3,-2){\circle*{0.1}}
\put(3,-2.8){\vdots}
\qbezier(3,1)(4.8,3)(5,1)
\qbezier(3,1)(4.8,-1)(5,1)
\put(3.2,1){$v$}
\put(3.2,0){$w_1$}
\put(3.2,-1){$w_2$}
\put(3.2,-2){$w_3$}
\put(1,0.2){$e_1$}
\put(5.1,1){$e_2$}
\put(0,-0.1){$>$}
\put(4.2,1.9){$>$}
\end{picture}}
\vspace{3 cm}

In this ultragraph, $\GG^1=\{e_1, e_2\}$, $s(e_1)=u$, $r(e_1)=\{v,w_1,w_2...\}$, $v=s(e_2)=r(e_2)$, and each $w_i$ is a sink. In this case, $\pinf$ contains two elements, $e_1e_2e_2...$ and $e_2e_2...$, and $\pfin=\{(e_1,w_i):i\in \N\}\cup\{(w_i,w_i):i\in \N\}$.
\end{exemplo}

\begin{definicao}
For an element $(\alpha,v)\in \mathfrak{p}^*$ we define the range and source maps by $r(\alpha,v)=v$ and $s(\alpha,v)=s(\alpha)$. In particular, for a sink $v$, $s(v,v)=v=r(v,v)$. We also extend the length map to the elements $(\alpha,v)$ by defining $|(\alpha,v)|:=|\alpha|$.
\end{definicao}

\section{A model for permutative, irreducible representations of $L_R(\mathcal{G})$ }

In this section, motivated by result in \cite{chen} for Leavitt path algebras, we define an irreducible representation associated to any ultragraph algebra $L_R(\mathcal{G})$. As we will see later using branching system theory, this representations models permutative, perfect, and irreducible representations of $L_R(\mathcal{G})$.

Recall that, unless stated otherwise, $R$ is a commutative unital ring.

\begin{definicao}
Two elements $\alpha,\beta\in \mathfrak{p}^*\cup \mathfrak{p}^\infty$ are equivalent if:
\begin{enumerate}
    \item $\alpha,\beta\in \mathfrak{p}^\infty$ and there are $i,j$ such that $\alpha_{i+k}=\beta_{j+k}$ for each $k\in \N$, or
    \item $\alpha,\beta \in \pfin$, where $\alpha=(a,v)$ and $ \beta=(b,v)$.
\end{enumerate}
\end{definicao}

\begin{remark}
For $\alpha\in\pfin\cup \pinf$, we denote by $[\alpha]$ the set of all the paths equivalent to $\alpha$, and by $\widetilde{\pinf}$ and $\widetilde{\pfin}$ respectively the set of equivalent classes of $\pinf$ and $\pfin$. Notice that $\widetilde{\pfin}\cap \widetilde{\pinf}=\emptyset$. Moreover, each class in $\pfin$ is given by a vertex which is a sink, so that the cardinality of $\widetilde{\pfin}$ and $G_s^0$ is the same.
\end{remark}

\begin{exemplo}\label{ex2}
In the following ultragraph

\vspace{2cm}

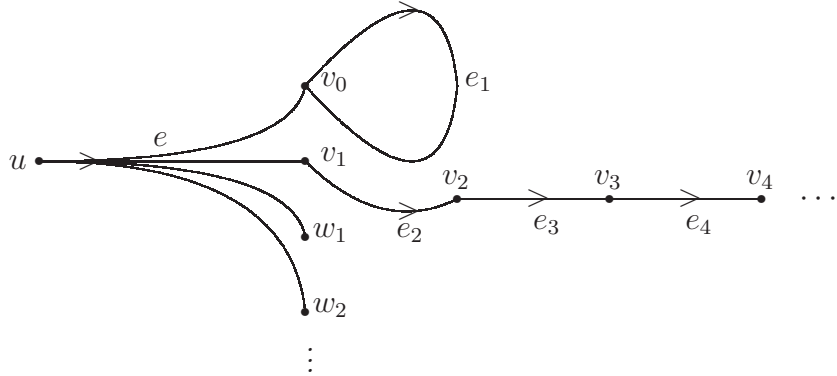
\centerline{
\setlength{\unitlength}{1cm}
\begin{picture}(9,0)
\put(-0.5,0){\circle*{0.1}}
\put(-0.9,-0.1){$u$}
\qbezier(-0.5,0)(2.8,0)(3,1)
\qbezier(-0.5,0)(2.8,0)(3,0)
\qbezier(-0.5,0)(2.8,0)(3,-1)
\qbezier(-0.5,0)(2.8,0)(3,-2)
\put(3,1){\circle*{0.1}}
\put(3,0){\circle*{0.1}}
\put(3,-1){\circle*{0.1}}
\put(3,-2){\circle*{0.1}}
\put(3,-2.8){\vdots}
\qbezier(3,1)(4.8,3)(5,1)
\qbezier(3,1)(4.8,-1)(5,1)
\qbezier(3,0)(4,-1)(5,-0.5)
\put(5,-0.5){\circle*{0.1}}
\qbezier(5,-0.5)(6,-0.5)(9,-0.5)
\put(7,-0.5){\circle*{0.1}}
\put(9,-0.5){\circle*{0.1}}
\put(4.8,-0.3){$v_2$}
\put(6.8,-0.3){$v_3$}
\put(8.8,-0.3){$v_4$}
\put(3.2,1){$v_0$}
\put(3.2,0){$v_1$}
\put(3.1,-1){$w_1$}
\put(3.1,-2){$w_2$}
\put(1,0.2){$e$}
\put(5.1,1){$e_1$}
\put(0,-0.1){$>$}
\put(4.2,1.9){$>$}
\put(4.2,-0.76){$>$}
\put(4.2,-1){$e_2$}
\put(5.9,-0.6){$>$}
\put(6,-0.9){$e_3$}
\put(7.9,-0.6){$>$}
\put(8,-0.9){$e_4$}
\put(9.5,-0.6){$\cdots$}
\end{picture}}
\vspace{3 cm}
we have, for example, that $ee_2e_3e_4...$ and $e_7e_8e_9...$ are equivalent, and so are the elements $(e,w_i)$ and $(w_i,w_i)$ for each $i$. There are only two classes in $\widetilde{\pinf}$, the class of $e_1e_1e_1...$ and the class of $ee_2e_3e_4...$. The set $\widetilde{\pfin}$ contains infinitely many elements, more specifically, $\widetilde{\pfin}=\{[(e,w_i)]:i\in \N\}$.
\end{exemplo}

\begin{definicao}\label{freemodule}
Let $\GG$ be an ultragraph. We denote by $\mathbb{P}$ be the free $R$-module generated by the basis $\{b_\alpha:\alpha\in \pfin\cup \pinf\}$. 
For an element $\alpha\in \pfin\cup \pinf$, we let  $\mathbb{P}_{\alpha}$ denote the submodule of $\mathbb{P}$ generated by $b_\alpha$, and $\mathbb{P}_{[\alpha]}$ denote the submodule of $\mathbb{P}$ generated by the elements $b_\beta$ with $\beta\in [\alpha]$. 
\end{definicao}

Notice that $\mathbb{P}_{[\alpha]}=\bigoplus\limits_{\beta\in [\alpha]}\mathbb{P}_\beta$, and $  \displaystyle \mathbb{P}=\bigoplus_{[\alpha]\in \widetilde{\pinf}}\mathbb{P}_{[\alpha]}\bigoplus\limits_{[\beta]\in \widetilde{\pfin}}\mathbb{P}_{[\beta]}.$

\begin{exemplo} For the ultragraph of Example \ref{ex1} notice that $\widetilde{\pinf}=\{[e_1e_2e_2...]\}$. Since $[e_1e_2e_2...]$ contains exactly two elements, then $\mathbb{P}_{[e_1e_2e_2...]}=R^2$. Moreover, for each $[(e_1,w_i)]\in \widetilde{\pfin}$, since $[(e,w_i)]$ also contains exactly two elements, we have that $\mathbb{P}_{[(e_1,w_i)]}=R^2$. Then

$$\mathbb{P}=\bigoplus_{[\alpha]\in \widetilde{\pinf}}\mathbb{P}_{[\alpha]}\bigoplus\limits_{[\beta]\in \widetilde{\pfin}}\mathbb{P}_{[\beta]}=R^2\bigoplus_{i\in \N}R^2=R^\N.$$

\end{exemplo}

Our aim is to define a representation $\pi:L_R(\GG)\rightarrow End_R(\mathbb{P})$, where $End_R(\mathbb{P})$ is the set of all $R$-endomorphism on $\mathbb{P}$, which is an  $R$-algebra. With this in mind, let us define some special elements in $End_R(\mathbb{P})$ as follows:
\begin{enumerate}
    \item for each element $A\in \GG^0$ define $P_A:\mathbb{P}\rightarrow \mathbb{P}$ by $P_A(b_{\alpha})=[s(\alpha)\in A]b_{\alpha}$. 
\item for each $e\in \GG^1$ define $S_e:\mathbb{P}\rightarrow \mathbb{P}$ by  $S_e(b_{\alpha})=[s(\alpha)\in r(e)]b_{e\alpha}$. 
\item for each $e^*\in(\GG^1)^*$, define $S_e^*:\mathbb{P}\rightarrow \mathbb{P}$ $S_e^*(b_{\alpha})=[\alpha_1=e]b_{\alpha_2\alpha_3...}$.

\end{enumerate}

\begin{remark} In the previous definition, the notation $[q]$ means $[q]=1$ if the statement $q$ is true and $[q]=0$ otherwise. In the second item, if $\alpha=(a,v)\in \pfin$ with $s(\alpha)\in r(e)$ then $e\alpha:=(ea,v)$. Particularly, if $\alpha=(v,v)$ and $v\in s(e)$ then $e\alpha=(e,v)$. Finally, in the third item, if $\alpha=(a,v)$ with $a_1=e$ then $S_e^*(b_\alpha)=b_{(a_2...a_{|a|},v)}$. In particular, if $\alpha=(e,v)$ then $S_e^*(b_\alpha)=b_{(v,v)}$, and if $\alpha=(v,v)$ then $S_e^*(b_\alpha)=0$ .
\end{remark}

The above endomorphisms induce a representation of $L_R(\GG)$ as described below.

\begin{teorema}\label{rep1} Let $\GG$ be an ultragraph.
There exists a representation $\pi:L_R(\GG)\rightarrow End_R(\mathbb{P})$ such that $\pi(p_A)=P_A$ for all $A\in \GG^0$, $\pi(s_e)=S_e$ for each edge $e$, and $\pi(s_e^*)=S_e^*$ for each $e^*\in (\GG^1)^*$.
\end{teorema}

\demo The proof of this theorem follows from the universality of $L_R(\GG)$ and from the definitions of $P_A, S_e$ and $S_e^*$.
\fim

\begin{remark} The representation $\pi$ above is not always faithful. For example let $\GG$ be the ultragraph of Example \ref{ex2}. Then for each $x\in \pfin\cup \pinf$ we have that 

\noindent $\pi(p_{v_0})(b_x)=P_{v_0}(b_x)=[s(x)=v_0]b_x=[x=e_1e_1e_1...]b_x$,  and $$\pi(s_{e_1})(b_x)=S_{e_1}(b_x)=[s(x)\in r(e_1)]b_{e_1x}=[s(x)=v_0]b_{e_1x}=[x=e_1e_1e_1...]b_{e_1x},$$ so that $ \pi(p_{v_0})(b_x)=\pi(s_{e_1})(b_x).$  Therefore $\pi(p_{v_0})=\pi(s_{e_1})$, and so $\pi$ is not faithful.

We will show later, in Theorem \ref{faithfulrep}, a sufficient and necessary condition for faithfulness of the representation $\pi$ of Theorem \ref{rep1}.
\end{remark}

Notice that for each $[p]\in \widetilde{\pfin}\cup \widetilde{\pinf}$, it holds that $\pi(P_A)(\mathbb{P}_{[p]})\subseteq \mathbb{P}_{[p]}$, $\pi(S_e)(\mathbb{P}_{[p]})\subseteq \mathbb{P}_{[p]}$ and $\pi(S_e^*)(\mathbb{P}_{[p]})\subseteq \mathbb{P}_{[p]}$, for each edge $e$ and $A\in \GG^0$. Therefore $\pi(L_R(\GG))(\mathbb{P}_{[p]})\subseteq \mathbb{P}_{[p]}$, and then we may consider the restriction of $\pi$ to $\mathbb{P}_{[p]}$, which is a new representation (that we still denote by $\pi$).

\begin{proposicao}\label{irredonpaths} Let $R$ be a field. For each $[p]\in \widetilde{\pfin}\cup \widetilde{\pinf}$ the representation $\pi:L_R(\GG)\rightarrow End_R(\mathbb{P}_{[p]})$ is irreducible.
\end{proposicao}

\demo First suppose that $[p]\in \widetilde{\pinf}$. Let $x,z\in [p]$. Then there are paths $\alpha,\beta\in \GG^*$ such that $x=\alpha \xi$ and $z=\beta \xi$ with $\xi\in \pinf$ and so $\pi(s_\alpha)\pi(s_\beta^*)(b_z)=b_x$.

Let $0\neq Y\subseteq \mathbb{P}_{[p]}$ be an invariant subspace and let $0\neq y=\sum \lambda_ib_{x_i}\in Y$ with $\lambda_i\neq 0$ for each $i$ and $x_i\neq x_j$ for each $i\neq j$. Since all the $x_i$ are distinct, there exists distinct $\alpha_i$, all of same length, such that $x_i=\alpha_i x_i'$ for each $i$. Now, for a fixed $j$, we get $$\pi(s_{\alpha_j}^*)(y)=\sum\limits_i S_{\alpha_j^*}(\lambda_ib_{x_i})=\sum\limits_i \lambda_iS_{\alpha_j^*}(b_{\alpha_ix_i'})=\lambda_jb_{x_j'}.$$

Since $R$ is a field then $b_{x_j'}\in Y$. By the first paragraph of this proof, we get that $b_x\in Y$ for each $x\in [p]$.

Now let $[p]\in \pfin$, and let $0\neq Y\subseteq \mathbb{P}_{[p]}$ be an invariant subspace. Let $0\neq y=\sum\limits_i\lambda_i b_{x_i}$ with each $\lambda_i\neq 0$ and $x_i=(\alpha_i,v)$ for each $i$, with $\alpha_i\neq \alpha_j$ for $i\neq j$. Fix an element $\alpha_{j_0}$ such that $|\alpha_{j_0}|\geq |\alpha_i|$ for each $i$. Note that $\pi(s_{\alpha_{j_0}}^*)(y)=\lambda_{j_0}b_{(v,v)}$, and hence $b_{(v,v)}\in Y$. Since for each $(\alpha,v)\in [p]$ we have that $\pi(s_{\alpha})(b_{(v,v)})=b_{(\alpha,v)}$, we conclude that $\mathbb{P}_{[p]}\subseteq Y$.  \fim

Notice that for each $[p]\in \widetilde{\pfin}\cup\widetilde{\pinf}$, since $\mathbb{P}_{[p]}$ is $\pi$-invariant, we can endow $\mathbb{P}_{[p]}$ with a left $L_R(\GG)$ module structure, where the product is defined by $a.b:=\pi(a)(b)$, for each $a\in L_R(\GG)$ and $b\in \mathbb{P}_{[p]}$. So we can consider $End_{L_R(\GG)}(\mathbb{P}_{[p]})$, the $R$-module of all endomorphisms of the $L_R(\GG)$ module $\mathbb{P}_{[p]}$. We then have the following extension of the first items in Theorems 3.3 and 3.7 of \cite{chen} to ultragraph Leavitt path algebras.

\begin{proposicao}
For each $[p]\in \widetilde{\pfin}\cup\widetilde{\pinf}$, the $R$-module $End_{L_R(\GG)}(\mathbb{P}_{[p]})$ is isomorphic to $R$.
\end{proposicao}

\demo Suppose first that $[p]\in \widetilde{\pinf}$, and let $\varphi\in End_{L_R(\GG)}(\mathbb{P}_{[p]})$. Fix $q\in [p]$. Then $\varphi(b_q)=\sum\limits_{i=1}^n\lambda_i b_{q^i}$ with $\lambda_i\neq 0$ and $q^i\neq q^j$ for each $i,j$. Suppose that $q^j\neq q$ for some $j$. 
Chose an index $m$ such that $q_1^j...q_m^j\neq q_1...q_m$ and $q_1^j...q_m^j\neq q_1^i...q_m^i$ for each $i\neq j$. Then $$\lambda_j(b_{q^j_{m+1}q^j_{m+2}}...)=S_{q_1^j...q_m^j}^*(\varphi(b_q))=\varphi(S_{q_1^j...q_m^j}^*(b_q))=0,$$ which is impossible. Then $q^j=q$ for each $j$, and so $\varphi(b_q)=\lambda_q b_q$ for all $q\in [p]$.

Now, for $r\in [p]$, write $r=\alpha x$ and $p=\beta x$, where $\alpha$ and $\beta$ are finite paths. Then $\varphi(S_\alpha^*(b_r))=\varphi(b_x)=\lambda_x b_x$, and so $$\lambda_x b_r=\lambda_x S_\alpha(b_x)=S_\alpha(\varphi(S_\alpha^*(b_r)))=(\varphi(S_\alpha S_\alpha^*(b_r)))=\varphi(b_r).$$ Similarly, (using $S_\beta$ and $S_\beta^*$ instead of $S_\alpha$ and $S_\alpha^*$) we get $\varphi(b_p)=\lambda_x b_p$. This implies that there exists $\lambda_\varphi\in R$ such that  $\varphi(b)=\lambda_\varphi b$ for each $b\in \mathbb{P}_{[p]}$. 

Suppose next that $[p]\in \widetilde{\pfin}$. Let $\varphi\in End_{L_R(\GG)}$ and let $(\alpha,v)\in [p]$. Then $\varphi(b_{(\alpha,v)})=\sum\limits_{i=1}^n\lambda_ib_{(\alpha_i,v)}$, where $\alpha_i\neq \alpha_j$. Notice that

$\begin{array}{ll}
\varphi(b_{(\alpha,v)}) & =\varphi(S_\alpha S_{\alpha}^*(b_{(\alpha,v)}))=S_\alpha S_{\alpha}^*(\varphi(b_{(\alpha,v)}))\\ &=\sum\limits_{i=1}^n\lambda_iS_\alpha S_\alpha^*(b_{(\alpha_i,v)}) =\sum\limits_{i:\alpha_i=\alpha \beta_i}\lambda_ib_{(\alpha_i,v)}.
\end{array}
$

\noindent So we may suppose that $|\alpha_i|\geq |\alpha|$ for each $i$. If $|\alpha_j|>|\alpha|$ for some $j$, then $\varphi(S_{\alpha_j}^*(b_{(\alpha,v)}))= 0 $ and, on the other hand, 

$\begin{array}{ll} \varphi(S_{\alpha_j}^*(b_{(\alpha,v)})) & =S_{\alpha_j}^*(\varphi(b_{(\alpha,v)}))=\sum\limits_{i=1}^n\lambda_jS_{\alpha_j}^*(b_{(\alpha_i,v)})\\ & =\lambda_jb_{(v,v)}+\sum\limits_{i\neq j}\lambda_iS_{\alpha_j}^*(b_{(\alpha_i,v)})\neq 0.
\end{array}$

\noindent Therefore  $|\alpha_i|=|\alpha|$ for each $i$. If $\alpha_j\neq \alpha$ for some $j$ then $$0=S_{\alpha_j}^*(\varphi(b_{(\alpha,v)}))=\lambda_jb_{(v,v)}+\sum\limits_{i\neq j}\lambda_iS_{\alpha_j}^*(b_{(\alpha_i,v)})\neq 0.$$ Then $\alpha_i=\alpha$ for each $i$ and so $\varphi(b_{(\alpha,v)})=\lambda_\alpha b_{(\alpha,v)}$ for all $(\alpha,v)\in [p]$. 

Now, for each $(\beta,v)\in [p]$ it holds that $\lambda_\beta b_{(\beta,v)}=\varphi(b_{(\beta,v)})=S_\beta(\varphi(b_{(v,v)}))=S_\beta(\lambda_{(v,v)})b_{(v,v)}=\lambda_{(v,v)}b_{(\beta,v)}$, so that there exists $\lambda_\varphi\in R$ such that  $\varphi(b)=\lambda_\varphi b$ for each $b\in \mathbb{P}_{[p]}$.

We leave to the reader the verification that the map $End_{L_R(\GG)}(\mathbb{P}_{[p]})\ni \varphi \mapsto \lambda_\varphi \in R$ is an isomorphism. \fim

Next we will extend the remainder items of Theorems 3.3 and 3.7 in \cite{chen} to ultragraph Leavitt paht algebras. Before we proceed we recall the notion of equivalence between representations. 

\begin{definicao}\label{equivrepres} Let $\pi: L_R(\GG) \rightarrow Hom_K(M)$ and $\varphi: L_R(\GG) \rightarrow Hom_R(N)$ be representations of $L_R(\GG)$, where $M$ and $N$ are $R$-modules. We say that $\pi$ is equivalent to $\varphi$ if there exists an $R$-module isomorphism $U:M\rightarrow N$ such that the diagram 
\begin{displaymath}
\xymatrix{
M \ar[r]^{\pi(a)} \ar[d]_{U} &
M \ar[d]^{U} \\
N \ar[r]_{\varphi(a)} & N }
\end{displaymath}
commutes, for each $a\in L_R(\GG)$.
\end{definicao}

\begin{proposicao} For each $[p],[q]\in \widetilde{\pfin}\cup \widetilde{\pinf}$, with $[p]\neq [q]$, the representations $\pi:L_R(\GG)\rightarrow  End_R(\mathbb{P}_{[p]})$ and  $\pi:L_R(\GG)\rightarrow  End_R(\mathbb{P}_{[q]})$ are not equivalent.
\end{proposicao}

\demo Suppose that the representations are equivalent, that is, suppose there is an isomorphism $U:\mathbb{P}_{[p]}\rightarrow \mathbb{P}_{[q]}$ such that $\pi(a)\circ U=U\circ \pi(a)$ for each $a\in L_R(\GG)$. 

First we analyse the case $[p]\in \widetilde{\pinf}$ and $[q]\in \widetilde{\pfin}$. Let $x\in [p]$. Then $U(b_x)=\sum\limits_i\lambda_i b_{(\alpha_i,v)}$. Let $m>|\alpha_i| $ for each $i$ and let $\beta=x_1...x_m$. Then $$U(b_{x_{m+1}x_{m+2}...})=U(S_\beta^*b_x)=U(\pi(s_\beta^*)(b_x))=\pi(s_\beta^*)(U(b_x))=S_\beta^*(\sum\limits_i\lambda_ib_{(\alpha_i,v)})=0,$$ which is impossible since $U$ is injective.

Now let $[p]\in \widetilde{\pfin}$ and $[q]\in \widetilde{\pinf}$. Proceeding similarly to the previous case, we get a contradiction for $U^{-1}$.

Fix now $[p], [q]\in\widetilde{\pinf}$. Then $U(b_p)=\sum\limits_ib_{q^i}$. Since $[p]\neq [q]$ then there exists an $m$ such that $p_1...p_m\neq q_1^i...q_m^i$ for each $i$. Let $\beta=p_1...p_m$. Then $U(b_{p_{m+1}p_{m+2}...})=U(S_\beta^*(b_p))=S_\beta^*(U(b_p))=0$, which is impossible, since $U$ is injective.

Finally consider the case $[p],[q]\in \widetilde{\pfin}$. Let $(\alpha,u)\in [p]$. Then $U(b_{(\alpha,u)})=\sum\limits_ib_{(\beta_i,v)}$ and, since $[p]\neq [q]$, we have that $u\neq v$. Hence 

$\begin{array}{ll} U(b_{(u,u)})& =U(P_uS_\alpha^*(b_{(\alpha,u)}))=U(\pi(p_u)\pi(s_\alpha^*)(b_{(\alpha,u)}))\\
& =\pi(p_u)\pi(s_\alpha^*)(U(b_{(\alpha,u)}))=\sum\limits_i P_u(S_\alpha^*(b_{(\beta_i,v)}))=0,
\end{array}$

\noindent where the last equality follows from the fact that $u\neq v$ and, since $u$ is a sink, $P_u(b_{\alpha'})\neq 0 \text{ iff }  \alpha' = (u,u)$. So $U(b_{(v,v)})=0$,
which is impossible.

Therefore we have proved that if the representations $\pi:L_R(\GG)\rightarrow  End_R(\mathbb{P}_{[p]})$ and  $\pi:L_R(\GG)\rightarrow  End_R(\mathbb{P}_{[q]})$ are equivalent then $[p]=[q]$.
\fim

\section{Branching systems}\label{bs}

Iterated function systems and branching systems are widely used in the study of representations of algebras associated to combinatorial objects, see for example \cite{cgcdg, MR3941382,  MR3404559, FGKP, FGJKP, GLR1, GLR, GLR2, GR3, GR, MR2848777, MR2903145, GR4, dd1, HR}. Hence it is interesting to note that the representation $\pi$ of Theorem~\ref{rep1} can be constructed via branching systems. This point of view will allow us to apply results in the theory of branching systems  to characterize the representation $\pi$. For example, in this section we will completely characterize when $\pi$ is faithful, a result for ultragraphs that also improves Proposition~4.4 in \cite{chen} regarding Leavitt path algebras. Before we proceed we recall the following relevant definitions (as in \cite{dd1}).

\begin{definicao}\label{branchsystem}
Let $\mathcal{G}$ be an ultragraph, $X$ be a set and let $\{R_e,D_A\}_{e\in \mathcal{G}^1,A\in \mathcal{G}^0}$ be a family of subsets of $X$. Suppose that
\begin{enumerate}
\item $R_e\cap R_f =\emptyset$, if $e \neq f \in \mathcal{G}^1$;
\item $D_\emptyset=\emptyset, \ D_A \cap D_B= D_{A \cap B}, \text{ and } D_A \cup D_B= D_{A \cup B}$ for all $A, B \in \mathcal{G}^0$;
\item $R_e\subseteq D_{s(e)}$ for all $e\in \mathcal{G}^1$;
\item\label{D_v=cup_{e in s^{-1}(v)}R_e} $D_v= \bigcup\limits_{e \in s^{-1}(v)}R_e$, if $0 <\vert s^{-1}(v) \vert<\infty$; and
\item for each $e\in \mathcal{G}^1$, there exist two bijective maps, $f_e:D_{r(e)}\rightarrow R_e$ and $f_e^{-1}:R_e \rightarrow D_{r(e)}$, such that $f_e\circ f_e^{-1}=Id_{R_e}$ and $f_e^{-1}\circ f_e=Id_{D_{r(e)}}$.
\end{enumerate}

We call $\{D_A,R_e,f_e\}_{e \in \mathcal{G}^1,A \in \mathcal{G}^0}$ a $\mathcal{G}$-algebraic branching system on $X$ or, shortly, a $\GG$-branching system, and use the notation $X=\{D_A,R_e,f_e\}_{e\in \GG^1, A\in \GG^0}$ to denote this branching system.
\end{definicao}

Let $M(X)$ be the $R-$module of all maps from $X$ to $R$ with finite support. In Proposition~4.5 of \cite{dd1} it is shown that a branching system induces a representation of $L_R(\GG)$ in $End(M(X))$. For our purposes, we will consider the following branching system on $\pfin\cup \pinf$.

\begin{itemize}
    \item For each $A\in \GG^0$ let $B_A=\{x\in \pfin\cup \pinf: s(x)\in A\}$,
    \item For each edge $e$ define $L_e=\{x\in \pfin\cup\pinf:x_1=e\}$,
    \item For each edge $e$ define $f_e:B_{r(e)}\rightarrow L_e$ by $f_e(x)=ex.$
\end{itemize}

\begin{remark}\label{onebranchsys} In the previous definition, if $x=(\alpha,v)$ then $f_e(x)=(e\alpha, v)$, and if $x=(v,v)$ then $f_e(x)=(e,v)$. It is easy to see that $f_e$ is bijective, and that
$\{B_A, L_e, f_e\}_{e\in \GG^1, A\in \GG^0}$ is a $\GG$-branching system. 
\end{remark}

\begin{remark}
The above branching system can also be seen as a partial action of the free group on the edges of the ultragraph and be used to realized $L_R(\GG)$ as a partial skew group ring, see \cite{cgcdg, dd2, dd3}.
\end{remark}

Next we make precise the representation $\varphi$ of $L_R(\GG)$ induced by the branching system defined above on $\pfin\cup \pinf$.

\begin{proposicao}\label{rep2} There is a representation $\varphi:L_R(\GG)\rightarrow End_R(M(\pfin\cup\pinf))$ such that:
$\varphi(p_A)(\phi)=\chi_{B_A}.\phi$ (where $\chi_{B_A}$ is the characteristic map on $B_A$), $\varphi(s_e)(\phi)=\chi_{L_e}.(\phi\circ f_e^{-1})$,
and $\varphi(s_e^*)(\phi)=\chi_{B_{r(e)}}.(\phi\circ f_e)$.
\end{proposicao}

\demo First let $N(\pfin\cup\pinf)$ be the $R$-module of all the maps from $\pfin\cup\pinf$ to $R$. Clearly $M(\pfin\cup\pinf)$ is a sub-module of $N(\pfin\cup\pinf)$. From Proposition 4.5 of \cite{dd1} we get a representation $\pi:L_R(\GG)\rightarrow End_R(N(\pfin\cup\pinf))$ such that  $\varphi(p_A)(\phi)=\chi_{B_A}.\phi$ for each $A\in \GG^0$,  $\varphi(s_e)(\phi)=\chi_{L_e}.(\phi\circ f_e^{-1})$
and $\varphi(s_e^*)(\phi)=\chi_{B_{r(e)}}.(\phi\circ f_e)$ for each edge $e$ and $\phi\in N(\pfin\cup\pinf)$. Now, it is easy to see that $M(\pfin\cup \pinf)$ is $\pi$-invariant, and so we get the desired representation.
\fim

As we mentioned before, the representations $\varphi$ and $\pi$ of $L_R(\GG)$ are equivalent, a fact we prove after setting up a basis for $M(\pfin\cup\pinf)$ below.

For each $x\in \pfin\cup\pinf$ define $\delta_x \in M(\pfin\cup \pinf)$ by $\delta_x(y)=1$ if $y=x$ and $\delta_x(y)=0$ if $y\neq x$. Notice that $\{\delta_x:x\in \pfin\cup\pinf\}$ is a basis of $M(\pfin\cup\pinf)$.

\begin{proposicao}\label{equivrep} The representations $\pi:L_R(\GG)\rightarrow End_R(\mathbb{P})$ of Theorem \ref{rep1} and $\varphi:L_R(\GG)\rightarrow End_R(M(\pfin\cup\pinf))$ of Proposition~\ref{rep2} are equivalent.
\end{proposicao}

\demo Define the isomorphism $U:\mathbb{P}\rightarrow M(\pfin\cup\pinf)$ by $U(\sum\limits_i\lambda_ib_{x_i})=\sum\limits_i \lambda_i \delta_{x_i}$. To show that $\pi(a)=U^{-1}\circ\varphi(a)\circ U$ for each $a\in L_R(\GG)$ it is enough to show that $\pi(p_A)=U^{-1}\circ\varphi(p_A)\circ U$ for each $A\in \GG^0$, $\pi(s_e)=U^{-1}\circ\varphi(s_e)\circ U$ and $\pi(s_e^*)=U^{-1}\circ\varphi(s_e^*)\circ U$ for each $e\in \GG^1$.

Let $A\in \GG^0$ and $x\in \pfin\cup \pinf$. Then

$\begin{array}{ll} U^{-1}(\varphi(p_A)(U(b_x)))&=U^{-1}(\varphi(p_A)(\delta_x))=U^{-1}(\chi_{B_A}(\delta_x))\\& =[x\in B_A]U^{-1}(\delta_x)=[x\in B_A]b_x=[s(x)\in A]b_x\\ & =\pi(p_A)(b_x),\end{array}$

\noindent and hence $\pi(p_A)=U^{-1}\circ\varphi(p_A)\circ U$.

Now let $e\in \GG^1$ and $x\in \pfin\cup\pinf$. Then 

$\begin{array}{ll} U^{-1}(\varphi(s_e)(U(b_x)))& =U^{-1}(\varphi(s_e)(\delta_x))=U^{-1}(\chi_{L_e}.(\delta_x\circ f_e^{-1})) \\ & =U^{-1}(\chi_{L_e}\delta_{f_e(x)})=[x\in B_{r(e)}]U^{-1}(\delta_{ex})\\ &= [x\in B_{r(e)}]b_{ex}=[s(x)\in r(e)]b_{ex}=\pi(s_e)(b_x),
\end{array}$

\noindent and therefore  $U^{-1}\circ\varphi(s_e)\circ U=\pi(s_e)$.

Analogously to what is done above one shows that $U^{-1}\circ\varphi(s_e^*) \circ U=\pi(s_e^*)$.

\fim

We are now ready to completely characterize faithfulness of the representation $\pi$ in terms of combinatorial properties of the underlying ultragraph, but first we recall the following definitions. 

\begin{definicao}[\cite{Tomsimplicity}] Let $\GG$ be an ultragraph. A closed path is a path
$\alpha \in \GG^*$ with $| \alpha | \geq 1$ and $s(\alpha)
\in r(\alpha)$. A closed path $\alpha$ is a cycle if $s(\alpha_i)\neq s(\alpha_j)$ for each $i\neq j$. An exit for a closed path is either an edge $e \in
\GG^1$ such that there exists an $i$ for which $s(e) \in
r(\alpha_{i})$ but $e \neq \alpha_{i+1}$, or a sink $w$ such that $w \in r(\alpha_i)$ for some $i$. We say that the ultragraph $\GG$ satisfies Condition~(L) if every closed path in $\GG$ has an
exit.
\end{definicao}

In Proposition 4.4 of \cite{chen}, the author shows that, for a row-finite graph $E$, Condition~$(L)$ is sufficient for faithfulness of the representation $\pi$ of $L_R(E)$. We show in the next theorem that the row-finite assumption is not necessary to describe faithfulness of $\pi$. In fact, using branching system theory, we show that for any ultragraph Condition~$(L)$ is necessary and sufficient for faithfulness of $\pi$.



\begin{teorema}\label{faithfulrep} Let $\GG$ be an ultragraph. Then the representation $\pi:L_R(\GG)\rightarrow End_R(\mathbb{P})$ of Theorem \ref{rep1} is faithful if, and only if, $\GG$ satisfies Condition~(L).
\end{teorema}

\demo By Proposition \ref{equivrep} it is enough to show that the representation $\varphi:L_R(\GG)\rightarrow End_R(M(\pfin\cup\pinf))$ is faithful if, and only if, $\GG$ satisfies Condition~$(L)$. 

Suppose that $\GG$ satisfies Condition~$(L)$. By Theorem 5.1 in \cite{dd1} $\varphi$ is faithful.

Now suppose that $\GG$ does not satisfy Condition~$(L)$. Then there is a closed path $c$ without exit. Notice that $B_{r(c)}=\{x\}$, where $x$ is the infinite path $x=ccc...$ and $f_c^n(x)=x$ for each $n\in\ \N$. Then, by Theorem 5.1 in \cite{dd1}, $\varphi$ is not faithful. \fim

\section{Perfect branching systems}

In this section we focus on perfect branching systems. Intuitively speaking a perfect branching system is one such that the Cuntz-Krieger relations, translated to the sets that form the branching system, holds for each non-sink vertex, and such that the whole set $X$ is the union of the "projection" sets associated to the vertices. This type of branching systems arise naturally, as the constructions in \cite{GLR1, GLR2, GR3, GR, MR2903145, GR4} show.

Our main goals in this section are to extend Lemma~5.4 and Theorem~5.6 in \cite{chen} to ultragraph Leavitt path algebras. In particular we show that there is always a morphism between a given branching system associated to an ultragraph and the branching system on $\pfin\cup\pinf$ described in the previous section. Furthermore, we will use this last result to characterize irreducible representations arising from branching systems. We start with the definition of a morphism between branching systems.

\begin{definicao}\label{isobranchsys} Let $\GG$ be an ultragraph. Let $X=\{D_A,R_e,f_e\}_{A\in \GG^0, e\in \GG^1}$ and $Y=\{B_A,L_e,g_e\}_{A\in \GG^0, e\in \GG^1}$ be two branching systems. A morphism from $X$ to $Y$ is a map $T:X\rightarrow Y$ such that $T(R_e)\subseteq L_e$ for each $e\in \GG^1$, $T(D_A)\subseteq B_A$ for each $A\in \GG^0$, and such that the diagram \begin{displaymath}
\xymatrix{
D_{r(e)} \ar[r]^{f_e} \ar[d]_{T} &
R_e \ar[d]^{T} \\
B_{r(e)} \ar[r]_{g_e} & L_e }
\end{displaymath}
commutes for each $e\in \GG^1$.

 The branchyng systems are isomorphic if there are mutually inverse morphisms $T:X\rightarrow Y$ and $T^{-1}:Y\rightarrow X$.

\end{definicao}

\begin{remark} From the previous definition we get that if $T:X\rightarrow Y$ is a morphism of branching systems then $T\circ f_e=g_e\circ T$. Now, composing this equality on the right with $f_e^{-1}$, and on the left with $g_e^{-1}$, we get $g_e^{-1}\circ T=T\circ f_e^{-1}$, so that the diagram \begin{displaymath}
\xymatrix{
R_e \ar[r]^{f_e^{-1}} \ar[d]_{T} &
D_{r(e)} \ar[d]^{T} \\
L_e \ar[r]_{g_e^{-1}} & B_{r(e)}}
\end{displaymath}
also commutes, for each $e\in \GG^1$.
\end{remark}


Let $\GG$ be an ultragraph and $X$ be a branching sytem. We denote by $N(X)$ the $R$-module of all the maps from $X$ to $R$. For two branching systems $X$ and $Y$ as in Definition~\ref{isobranchsys}, let $\pi:L_R(\GG)\rightarrow End_R(N(X))$ and $\varphi:L_R(\GG)\rightarrow End_R(N(Y))$ be the representation induced by these branching systems, as in Proposition 4.5 of \cite{dd1}. Recall that $\pi(p_A)(\phi)=\chi_{D_A}\phi$, $\pi(s_e)(\phi)=\chi_{R_e}(\phi\circ f_e^{-1})$, $\pi(s_e^*)(\phi)=\chi_{D_{r(e)}}(\phi\circ f_e)$, for each $\phi\in N(X)$, and analogous description holds for $\varphi$. 

Next we notice that isomorphic branching systems induce equivalent representations of $L_R(\GG)$, a result that follows directly from the lemma below. 

\begin{lema}\label{morphismlemma} Let $\GG$ be an ultragraph, $X=\{D_A,R_e,f_e\}_{A\in \GG^0, e\in \GG^1}$ and $Y=\{B_A,L_e,g_e\}_{A\in \GG^0, e\in \GG^1}$ be two branching systems, $T:X\rightarrow Y$ be a morphism of branching systems, and let $\pi:L_R(\GG)\rightarrow End_R(N(X))$ and $\varphi:L_R(\GG)\rightarrow End_R(N(Y))$ be the induced representations as described above. Suppose that $T^{-1}(L_e)=R_e$ for each edge $e$ and $T^{-1}(B_{A})=D_{A}$ for each $A\in \GG^0$, and let $U:N(Y)\rightarrow N(X)$ be defined by $U(\phi)=\phi\circ T$. Then $\pi(a)\circ U=U\circ \varphi(a)$ for each $a\in L_R(\GG)$.
\end{lema}

\demo To show that $\pi(a)\circ U=U\circ \varphi(a)$ for each $a\in L_R(\GG)$ it is enough to verify this equality for $a=s_e$, $a=s_e^*$ and $a=p_A$, for each edge $e$ and $A\in \GG^0$.
Let $e\in \GG^1$. For each $\phi\in N(X)$, $$\pi(s_e)(U(\phi)))=\pi(s_e)(\phi\circ T)=\chi_{R_e}.(\phi\circ T\circ f_e^{-1})=\chi_{T^{-1}(L_e)}.(\phi\circ g_e^{-1}\circ T)=$$
$$=\chi_{(L_e)}\circ T.(\phi\circ g_e^{-1}\circ T)=(\chi_{L_e}.\phi\circ g_e^{-1})\circ T=U(\varphi(\phi)),$$
 and so $\pi(s_e)\circ U=U\circ \varphi(s_e)$.

Analogously one shows that $\pi(s_e^*)\circ U=U\circ \varphi(s_e^*)$ and $\pi(p_A)\circ U=U\circ \varphi(p_A)$ for each edge $e$ and $A\in \GG^0$.
\fim

\begin{corolario}Let $\GG$ be an ultragraph, $X,Y$ be two isomorphic branching systems, and let $\pi:L_R(\GG)\rightarrow End_R(N(X))$ and $\varphi:L_R(\GG)\rightarrow End_R(N(Y))$ be their induced representations. Then $\pi$ and $\varphi$ are equivalent. 
\end{corolario}

Let $X$ and $Y$ be two isomorphic branching systems of an ultragraph $\GG$, with the branching system isomorphism $T:X\rightarrow Y$. Let $U:N(Y)\rightarrow N(X)$ be the induced isomorphism of $R$-modules, defined by $U(\phi)=\phi\circ T$ (recall that $N(Y)$ is the set of all the maps from $Y$ to $R$). Notice that $M(Y)$ (the set of all the maps from $Y$ to $R$ with finite support) is isomorphic to $M(X)$ via $U$. Moreover, $M(X)$ is $\pi$-invariant where $\pi$ is as in the previous corollary, so that we may consider the restricted representation $\pi:L_R(\GG)\rightarrow End_R(M(X))$ and similarly we get $\varphi:L_R(\GG)\rightarrow End_R(M(Y))$. By combining those facts, we get the following:

\begin{corolario}\label{equivrepbs} Let $\GG$ be an ultragraph, $X,Y$ be two isomorphic branching systems, and let $\pi:L_R(\GG)\rightarrow End_R(M(X))$ and $\varphi:L_R(\GG)\rightarrow End_R(M(Y))$ be their induced representations. Then $\pi$ and $\varphi$ are equivalent. 
\end{corolario}

As we mentioned before our aim in this section is to study perfect branching systems. We make this definition precise below (notice that this generalizes the definition given in Section~5 of \cite{chen}).

\begin{definicao} Let $\GG$ be an ultragraph, and $X$ be a $\GG-$algebraic branching system. We say that $X$ is perfect if $X=\bigcup\limits_{v\in G^0}D_v$, and $X_v=\bigcup\limits_{e\in s^{-1}(v)}X_e$ for each non-sink $v\in G^0$.  
\end{definicao}

\begin{exemplo}\label{perfectbs} The $\GG$-branching system on $\pfin\cup \pinf$ of Section~\ref{bs}, namely $\{B_A,L_e,f_e\}_{A\in \GG^0, e\in \GG^1}$, where $B_A=\{x\in \pfin\cup \pinf: s(x)\in A\}$, $L_e=\{x\in \pfin\cup\pinf:x_1=e\}$ and $f_e:B_{r(e)} \rightarrow L_e$ defined by $f_e(x)=ex$ is perfect.
\end{exemplo}

\begin{remark} If $X=\{D_A,R_e,f_e\}_{e\in \GG^1, A\in \GG^0}$ is a perfect branching system then $D_A=\bigcup\limits_{v\in A}D_v$. To see this,  first note that, for each $A\in \GG^0$, we have $D_v\cap D_A=D_v$ for each $v\in A$, so that $\bigcup\limits_{v\in A}D_v\subseteq D_A$. Moreover, if $x\in D_A$ then $x\in D_u$ for some $u$, since $X=\bigcup\limits_{v \in G^0}D_v$. If we suppose that that $u\notin A$ then we get that $x\in D_u\cap D_A=\emptyset$, which is impossible. Therefore $D_A=\bigcup\limits_{v\in A}D_v$.
\end{remark}

Morphisms from perfect branching systems have a special property, which we record below.


\begin{lema}\label{lemmaperfect}
Let $X=\{D_A, R_e, f_e\}_{e\in \GG^1, A\in \GG^0}$ and $Y=\{B_A, L_e, g_e\}_{e\in \GG^1, A\in \GG^0}$ be branching systems of an ultragraph $\GG$, let $T:X\rightarrow Y$ be a morphism, and suppose that $X$ is perfect. Then $T^{-1}(L_e)=R_e$ and $T^{-1}(B_A)=D_A$, for each edge $e$ and $A\in \GG^0$. 
\end{lema}

\demo To see that $T^{-1}(L_e)=R_e$, first note that, since $T(R_e)\subseteq L_e$ then $R_e\subseteq T^{-1}(L_e)$. Moreover, if $x\in T^{-1}(L_e)\setminus R_e$ then, since $X$ is perfect, $x\in D_v$ for some vertex $v$. If $s(e)=v$ then, since $X$ is perfect, $x\in R_{e_0}$ for some $e_0\neq e$, and so $T(x)\in L_{e_0}$. Since $L_{e_0}\cap L_e=\emptyset$ then $T(x)\notin L_e$, which is impossible since $x\in T^{-1}(L_e)$. If $s(e)\neq v$ then, since $T(x)\in B_v$ and $B_v\cap B_{s(e)}=\emptyset$, we have $T(x)\notin L_e$, which is also impossible. Therefore $T^{-1}(L_e)=R_e$. 

To verify that $T^{-1}(B_A)=D_A$, note first that $D_A\subseteq T^{-1}(B_A)$. If $x\in T^{-1}(B_A)\setminus D_A$, then $x\in D_v$ for some vertex $v$ (since $X$ is perfect), and so $T(x)\in B_v$. If $v\notin A$ then $B_v\cap B_A=\emptyset$ and so $T(x)\notin B_A$, which is impossible. If $v\in A$ then, since $X$ is perfect, we have that $D_A=\bigcup\limits_{v\in A}D_v$ and hence $x\in D_A$, which is also impossible. Therefore $T^{-1}(B_A)=D_A$.
\fim

Joining the lemma above and Lemma~\ref{morphismlemma} we get the following.

\begin{proposicao}\label{equivrep8}  Let $\GG$ be an ultragraph, $X,Y$ be two branching systems, $T:X\rightarrow Y$ be a morphism of branching systems, and let $\pi:L_R(\GG)\rightarrow End_R(N(X))$ and $\varphi:L_R(\GG)\rightarrow End_R(N(Y))$ be the induced representations. Suppose that $X$ is perfect and let $U:N(Y)\rightarrow N(X)$ be defined by $U(\phi)=\phi\circ T$. Then $\pi(a)\circ U=U\circ \varphi(a)$ for each $a\in L_R(\GG)$.
\end{proposicao}

\begin{remark}
Notice that the above proposition does not necessarily imply that the representations are equivalent, since the map $U$ may not be an isomorphism.
\end{remark}

We now describe a relationship between perfect branching systems and the branching system on  $\pfin\cup\pinf$ from Example~\ref{perfectbs}.

\begin{proposicao}\label{morphism} Let $\GG$ be an ultragraph and $X=\{R_e, D_A, g_e\}_{e\in \GG^1, A\in \GG^0}$ be a perfect $\GG$-branching system. Then there exists a morphism from $X$ to the branching system $\pfin\cup\pinf=\{L_e, B_A,f_e\}_{e\in \GG^1, A\in \GG^0}$ of Example~\ref{perfectbs}. 
\end{proposicao}

\demo First we define a map $T:X\rightarrow \pfin\cup \pinf$.

Let $x\in X$. Then $x\in D_{v_1}$ for some vertex $v_1$. If $v_1$ is a sink define $T(x)=(v_1,v_1)$. If $v_1$ is not a sink then $D_{v_1}=\bigcup\limits_{e\in s^{-1}(v_1)}R_e$. Let $e_1\in s^{-1}(v_1)$ be such that $x\in R_{e_1}$, and consider the element $g_{e_1}^{-1}(x)$. Notice that $g_{e_1}^{-1}(x)\in D_{v_2}$ for some vertex $v_2 \in r(e_1)$. If $v_2$ is a sink then define $T(x)=(e_1,v_2)$, otherwise there exists an edge $e_2$ such that $g_{e_1}^{-1}(x)\in R_{e_2}$. Consider the element $g_{e_2}^{-1}(g_{e_1}^{-1}(x))$, which belongs to $D_{v_3}$ for some vertex $v_3 \in r(e_2)$. If $v_3$ is a sink define $T(x)=(e_1e_2,v_3)$, otherwise there is an edge $e_3$ such that $g_{e_2}^{-1}(g_{e_1}^{-1}(x))\in R_{e_3}$. Proceeding recursively, we define $T(x)$ either as the element $(e_1e_2...e_n,v_{n+1})\in \pfin$ or the element $e_1e_2e_3...\in \pinf$. 

Notice that from the definition of $T$ we have $T(R_e)\subseteq L_e$ for each edge $e$. Moreover, if $v$ is a sink then $T(D_v)=\{(v,v)\}=B_v$, and if $v$ is not a sink then $T(D_v)=T(\bigcup_{e\in s^{-1}(v)}R_e)=\bigcup_{e\in s^{-1}(v)}T(R_e)\subseteq \bigcup_{e\in s^{-1}(v)}L_e=B_v$. Now, for $A\in \GG^0$, $T(D_A)=T(\bigcup\limits_{v\in A}D_v)=\bigcup\limits_{v\in A}T(D_v)\subseteq \bigcup\limits_{v\in A}B_v=B_A$.

Let $e$ be an edge and $x\in L_e$. Then, from the definition of $T$, we get that $T(g_e(x))=eT(x)=f_e(T(x))$, and so $T$ is a morphism of branching systems.

\fim



The morphism $T$ of the previous proposition is not always injective nor surjective. For example, let $\GG$ be the ultragraph with one edge $e$ and two vertices $u,v$, where $s(e)=u$ and $r(e)=\{u,v\}$.

\vspace{1cm}
\centerline{
\setlength{\unitlength}{0.5cm}
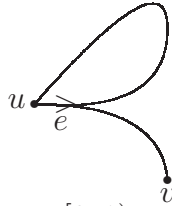
\begin{picture}(5,0)
\put(-0.5,0){\circle*{0.2}}
\put(-1.2,-0.1){$u$}
\put(2.8,-2.6){$v$}
\put(3,-2){\circle*{0.2}}
\put(0,-0.2){$>$}
\put(0,-0.7){$e$}
\qbezier(-0.5,0)(2.99,-0.3)(3,2)
\qbezier(-0.5,0)(2.99,4)(3,2)
\qbezier(-0.5,0)(2.8,0)(3,-2)
\end{picture}}
\vspace{1 cm}
Let $X=[0,2)$, define $R_e=[0,1)=D_u$, $D_v=[1,2)$ and $D_{r(e)}=[0,2)$, and let $f_e:D_{r(e)}\rightarrow R_e$ be any bijection. Then $X$ is a perfect branching system and for this ultragraph, $\pfin\cup\pinf=\{eee...,(e,v)\}$. Since $X$ is infinite and $\pfin\cup\pinf$ is a finite set, the morphism $T:X\rightarrow \pfin\cup\pinf$ is not injective.     

For an example where the morphism $T$ is not surjective, let $\GG$ be a graph with two loops $e_1$ and $e_2$ based on a vertex $u$.

\vspace{1cm}
\centerline{
\setlength{\unitlength}{1cm}
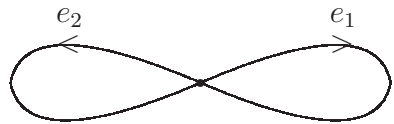
\begin{picture}(0,0)
\put(-0.5,0){\circle*{0.1}}
\qbezier(-0.5,0)(1.8,-1)(2,0)
\qbezier(-0.5,0)(1.8,1)(2,0)
\put(1.2,0.4){$>$}
\put(1.2,0.8){$e_1$}
\qbezier(-0.5,0)(-2.8,-1)(-3,0)
\qbezier(-0.5,0)(-2.8,1)(-3,0)
\put(-2.4,0.4){$<$}
\put(-2.4,0.8){$e_2$}
\end{picture}}
\vspace{1 cm}

Let $R_{e_1}$ and $R_{e_2}$ be two infinite countable disjoint sets,  $X=D_u=R_{e_1}\cup R_{e_2}$ and let $f_{e_i}:D_u\rightarrow R_{e_i}$ be a fixed bijection, for $i\in \{1,2\}$. Notice that this branching system is perfect. Moreover, $\pfin\cup \pinf=\pinf$ is not countable. Therefore the morphism $T$ of the previous proposition is not surjective. 

Although the morphism $T$ is not always surjective, we get the following lemma, which will be used in the next theorem.

\begin{lema}\label{surjectivelemma} Let $\GG$ be an ultragraph, $X=\{D_A, R_e, g_e\}_{e\in \GG^1, A\in \GG^0}$ be a perfect $\GG$-branching system, and let $T:X\rightarrow \pfin\cup\pinf$ be the morphism of Proposition \ref{morphism}. If $p\in \pfin\cup\pinf$ belongs to $T(X)$ then $[p]\subseteq T(X)$.
\end{lema}

\demo Since $T$ is a morphism we have that $T\circ g_e=f_e\circ T$ and $f_e^{-1}\circ T=T\circ g_e^{-1}$ for each edge $e$, and therefore it holds that $T\circ g_\alpha=f_\alpha\circ T$ and $f_\alpha^{-1}\circ T=T\circ g_\alpha^{-1}$, for each path $\alpha$. Now suppose $p\in \pinf \cap T(X)$, and write $p=T(x)$ for some $x\in X$. Let $y\in [p]$ and write $y=\alpha c$, where $p=\beta c$ and $c\in \pinf$. Then $c=f_\beta^{-1}(p)=f_\beta^{-1}(T(x))=T(g_\beta^{-1}(x))$ so that $c\in T(X)$. Let $c=T(d)$, where $d\in X$. Then $y=\alpha c=f_\alpha(c)=f_\alpha(T(d))=T(g_\alpha(d))$, and so $y\in T(X)$. Therefore $T(X)=[q]$. Similarly one shows that $T(X)=[q]$ if $q\in \pfin$.
\fim

\begin{remark}\label{restrictedbranchsys} Recall the branching system of Example~\ref{perfectbs}. Notice that for each $p\in \pfin\cup\pinf$ and $e\in \GG^1$ it holds that $f_e^{-1}(L_e\cap [p])\subseteq [p]$ and $f_e(B_{r(e)}\cap [p])\subseteq [p]$. Therefore we get a new branching system in $[p]$, by taking the intersections of $L_e$ and $B_A$ with $[p]$.
\end{remark}

We finish this section characterizing irreducible representations of $L_R(\GG)$ associated to perfect branching systems.

\begin{teorema}\label{irreducibletheorem} Let $\GG$ be an ultragraph, $R$ be a field, $X=\{D_A,R_e,g_e\}_{e\in \GG^1, A\in \GG^0}$ be a perfect $\GG$-branching system, and $\psi:L_R(\GG)\rightarrow End_R(M(X))$ be the associated representation. Then $\psi$ is irreducible if, and only if, $X$ is isomorphic to $[p]$ for some $[p]\in \widetilde{\pfin}\cup \widetilde{\pinf}$.
\end{teorema}

\demo If the branching system $X$ is isomorphic to the branching system $[p]$ for some $[p]\in \widetilde{\pfin}\cup  \widetilde{\pinf}$ (where $[p]$ is the branching system as in Remark \ref{restrictedbranchsys}) then $\psi$ is irreducible by Proposition~\ref{irredonpaths} and Corollary~\ref{equivrepbs}.

Suppose that $\psi$ is irreducible. Denote by $\varphi:L_R(\GG)\rightarrow End_R(M(\pfin\cup \pinf))$ the representation arising from the branching system on $ \pfin\cup \pinf$ defined on Example~\ref{perfectbs}. Let $T:X\rightarrow \pfin\cup\pinf$ be the morphism of branching systems defined in the proof of Proposition~\ref{morphism}. Notice that $T$ induces a map $V:M(X) \rightarrow M(\pfin\cup \pinf)$ that takes $\delta_x$ to $\delta_{T(x)}$, which is an $R$-homomorphism. We show that this map intertwines the representations, that is, $V \circ \psi (a) = \varphi(a) \circ V $ for all $a\in L_R(\GG) $. It is enough to verify that $(V \circ \psi (a))(\delta_x) = (\varphi(a) \circ V)(\delta_x)$ for each $x\in X$, and for $a=s_e$, $a=s_e^*$ and $a=p_A$ for each edge $e$ and $A\in \GG^0$.

Fix $A\in\GG^0$, and $x\in X$. Then $$V(\psi(p_A)(\delta_x))=[x\in D_A]V(\delta_x)= [x\in D_A]\delta_{T(x)},$$ and $$\varphi(p_A)(V(\delta_x))=\varphi(p_A)(\delta_{T(x)})=[T(x)\in B_A]\delta_{T(x)}.$$ Notice that $[T(x)\in B_A]=[x\in T^{-1}(B_A)]=[x\in D_A]$, where the last equality follows from Lemma \ref{lemmaperfect}. Therefore $V(\psi(p_A)(\delta_x))=\varphi(p_A)(V(\delta_x))$.

Next let $e\in \GG^1$ and $x\in X$. Then 

\begin{center}
$\begin{array}{ll} V(\psi(s_e)(\delta_x)) &= V(\chi_{R_e}.\delta_x\circ g_e^{-1})=[g_e(x)\in R_e]V(\delta_{g_e(x)})\\& =[g_e(x)\in R_e]\delta_{T(g_e(x))}=[x\in D_{r(e)}]\delta_{T(g_e(x))},
\end{array}$
\end{center}

\noindent and

\begin{center}
$ \begin{array}{ll}
\varphi(s_e)(V(\delta_x))& =\varphi(\delta_{T(x)})=\chi_{L_e}\delta_{T(x)}\circ f_e^{-1}\\ &=\chi_{L_e}(\delta_{f_e(T(x))})=[f_e(T(x))\in L_e]\delta_{f_e(T(x))}.
\end{array}$
\end{center}

Note that $f_e(T(x))=T(g_e(x))$ since $T$ is a morphism. Moreover 

\begin{center}
$     \begin{array}{ll}
[f_e(T(x))\in L_e]& =[T(x)\in f_e^{-1}(L_e)]=[T(x)\in B_{r(e)}]\\&=[x\in T^{-1}(B_{r(e)})]=[x\in D_{r(e)}],
\end{array}$
\end{center}
where the last equality follows from Lemma \ref{lemmaperfect}. 
So it follows that $V(\psi(s_e))(\delta_x)=\varphi(s_e)V(\delta_x)$ for every $x$ and hence $V\circ \psi(s_e)=\varphi(s_e)\circ V$ for all $e$.

Similarly to what is done above one shows that $V\circ \psi(s_e^*)=\varphi(s_e^*)\circ V$, and hence we conclude that $V\circ \psi(a)=\varphi(a)\circ V$ for each $a\in L_R(\GG)$.

Now, if $V$ is not injective, its kernel is invariant under $\psi$ (from the intertwining condition). Hence $V$ is injective and so is $T$. 

Finally, for each $p\in \pfin\cup \pinf$ let $Y_{[p]}\subseteq M(\pfin\cup \pinf)$ be the submodule generated by $\{\delta_x:x\in [p]\}$. It is easy to see that $Y_{[p]}$ is $\varphi$-invariant. Then, $V^{-1}(Y_{[p]})$ is $\psi$-invariant. Notice that $V(X)\cap Y_{[q]}\neq \emptyset$ for some $[q]$. Since $\psi$ is irreducible we have that $V^{-1}(Y_{[q]})=M(X)$, and so $V(M(X))\subseteq Y_{[q]}$. Therefore $T(X)\subseteq [q]$ and, from Lemma~\ref{surjectivelemma}, we get $T(X)=[q]$. \fim

\begin{remark}
In fact, in the above theorem, the assumption that $R$ is a field is only necessary to show the sufficient condition for irreducibility of $\psi$ (since we use Proposition~\ref{irredonpaths}).
\end{remark}

\section{Permutative representations}

Permutative representations of combinatorial algebras such as the Cuntz-Krieger, graph and ultragraph algebras have connections with the theory of operator algebras, dynamical systems, and pure algebra (see\cite{bratteli, kawamura8, GR, GLR1, FGKP}), and therefore are a subject of much interest. In this section we characterize the perfect, irreducible and permutative representations of an ultragraph Leavitt path algebra.

Let $\psi:L_R(\GG)\rightarrow End_R(M)$ be a representation, where $M$ is an $R-$module. Define the submodules $M_e=\psi(s_es_e^*)(M)$, for each edge $e$, and $M_A=\psi(p_A)(M)$ for each $A\in \GG^0$. Notice that:

\begin{enumerate}
    \item for each edge $e$, $\psi(s_e):M_{r(e)}\rightarrow M_e$ is invertible, with inverse $\psi(s_e^*)$;
    \item $M_u\cap M_v\{0\}$ and $M_e\cap M_f=0$ for each vertices $u\neq v$ and edges $e\neq f$;
    \item $M_v\supseteq \bigoplus\limits_{e\in s^{-1}(v)}M_e$ for each nonsink $v$ and if $0<|s^{-1}(v)|<\infty$ then $M_v= \bigoplus\limits_{e\in s^{-1}(v)}M_e$;
    \item $M_A\supseteq \bigoplus\limits_{v\in A}M_v$ for each $A\in \GG^0$.
    \end{enumerate}

\begin{definicao}
A representation $\psi:L_R(\GG)\rightarrow End_R(M)$ is called a perfect representation if $M_v= \bigoplus\limits_{e\in s^{-1}(v)}M_e$ for each nonsink $v$ and $M= \bigoplus\limits_{v\in G^0}M_v$.
\end{definicao}

From now on we suppose that $\psi$ is a perfect representation. 
Our goal is to construct a branching system associated to $\psi$. Below we describe how to define the sets of this branching system.

For each edge $e$ let $B_e$ be a basis of $M_e$. For each nonsink $v$, let $B_v=\bigcup\limits_{e\in s^{-1}(v)}B_e$. Since $\psi$ is perfect we have that $B_v$ is a basis of $M_v$. For each sink $v$, let $B_v$ be some basis of $M_v$.

Notice that $B=\bigcup\limits_{v\in \GG^0}B_v$ is a basis of $M$, since $\psi$ is perfect, and write $B=\{b_x: x\in X\}$, where $X$ is the index set of the basis $B$.

\begin{remark} From the hypothesis that $\psi$ is perfect it follows that $B_A=\bigcup\limits_{v\in A}B_v$ is a basis of $M_A$, for each $A\in \GG^0$.
\end{remark}

Next we define the subsets of $X$ that will form the desired branching system: For each edge $e$ write $B_e=\{b_x:x\in R_e\}$ where $R_e\subset X$ is the index set of the basis $B_e$, and for each $A\in \GG^0$ write $B_A=\{b_x:x\in D_A\}$ where $D_A\subset X$ is the index set of the basis $B_A$. Note that for edges $e\neq f$ we have $R_e\cap R_f=\emptyset$ (since $M_e\cap M_f=\{0\}$), and similarly $D_u\cap D_v=\emptyset$ for vertices $u\neq v$.

To define the bijections between the subsets defined above we need to restrict to permutative representations. We recall below Definition~6.1 of \cite{dd1}, already simplified to perfect representations.

\begin{definicao}\label{permcalor} Let $\psi:L_R(\GG)\rightarrow End_R(M)$ be a perfect representation. We say that $\psi$ is permutative if it is possible to choose basis $B_e$ and $B_v$ as described above and such that $\psi(s_e)(B_{r(e)})=B_e$.
\end{definicao}

\begin{remark} Notice that since $\psi(s_e):M_{r(e)}\rightarrow M_e$ is invertible, with inverse $\psi(s_e^*)$, we have that  $\psi(s_e)(B_{r(e)})=B_e$ is equivalent to $\psi(s_e^*)(B_e)=B_{r(e)}$. So $\psi$ is permutative if, and only if, for each edge $e$ it holds that $\psi(s_e^*)(B_e)=B_{r(e)}$.
\end{remark}

For a general permutative representation $\psi:L_R(\GG)\rightarrow End_R(M)$ we may define, for each edge $e$, the bijection $f_e:D_{r(e)}\rightarrow R_e$ by $f_e(x)=y$, where $\psi(s_e)(b_x)=b_y$. This leads us to the desired branching system, as we state below (and leave the proof to the reader).

\begin{proposicao}\label{bs1}
Let $\psi:L_R(\GG)\rightarrow End_R(M)$ be a perfect permutative representation. Then $\{R_e,D_A,f_e\}_{e\in \GG^1, A\in \GG^0}$ defined as above is a perfect branching system in $X$.
\end{proposicao}

A key example of a perfect permutative representation is the representation of Theorem~\ref{rep1}, as we see below.

\begin{exemplo}\label{perfectpermutativerep} The representation $\pi:L_R(\GG)\rightarrow End_R(\mathbb{P})$ obtained in Theorem \ref{rep1} is perfect and permutative. In fact, notice that in this case, for each edge $e$, $M_e$ is the submodule of $\mathbb{P}$ generated by $\{b_\alpha:\alpha\in \pfin\cup\pinf \text{ and }\alpha_1=e\}$, and for each $A\in \GG^0$, $M_A$ is the submodule of $\mathbb{P}$ generated by $\{b_\alpha:\alpha\in \pfin\cup\pinf \text{ and }s(\alpha)\in A\}$. It is easy to see that if $v$ is not a sink then $M_v=\bigoplus\limits_{e\in s^{-1}(v)}M_e$ and that $\mathbb{P}=\bigoplus\limits_{v\in G^0}M_v$, and hence $\pi$ is perfect. To see that $\pi$ is permutative, for each edge $e$, let $B_e\subseteq M_e$ be defined by $B_e=\{b_\alpha:\alpha\in \pfin\cup\pinf \text{ and }\alpha_1=e\}$ and let $B_{r(e)}\subseteq M_{r(e)}$ be the set $B_{r(e)}=\{b_\alpha:\alpha\in \pfin\cup\pinf \text{ and }s(\alpha)\in r(e)\}$. Now, it follows from the definition of $\pi(s_e)$ that $\pi(s_e):B_{r(e)}\rightarrow B_e$ is a bijection.
\end{exemplo}

Representations induced by the branching system of Proposition~\ref{bs1} form a model for perfect, permutative representations, as we show below. 

\begin{proposicao}\label{eqrep2} Let $\psi:L_R(\GG)\rightarrow End_R(M)$ be a perfect, permutative representation. Then $\psi$ is equivalent to the representation $\varphi: L_R(\GG)\rightarrow End_R(M(X))$ induced by the branching system of Proposition \ref{bs1}.
\end{proposicao}

\demo 
Let $B=\bigcup\limits_{v\in \GG^0}B_v$ be a basis of $M=\displaystyle \bigoplus_{v\in G^0} M_v$ as in the definition of a permutative representation (see Definition~\ref{permcalor}).
Following Theorem 6.5 in \cite{dd1}, it is enough to verify that $\psi(s_e^*)(b)=0$ for each edge $e$ and $b\in B\setminus B_e$, and $\psi(p_A)(b)=0$ for each $b\in B\setminus B_A$ for each $A\in \GG^0$. 

Fix an edge $e$ and $b\in B\setminus B_e$. Since $\psi$ is perfect, there exists a vertex $v\in G^0$ such that $b\in B_v$. If $v$ is a sink then $\psi(s_e^*)(b)=\psi(s_e^*)\psi(p_v)(b)=\psi(s_e^*)\psi(p_{s(e)})\psi(p_v)(b)=0$ since $s(e)\neq v$. If $v$ is not a sink then $B_v=\bigcup\limits_{f\in s^{-1}(v)}B_f$ and so $b\in B_f$, for some $f\neq e$. Then $b=\psi(s_f)\psi(s_f^*)(b)$ and so $\psi(s_e^*)(b)=\psi(s_e^*)\psi(s_fs_f^*)(b)=0$ since $\psi(s_e)^*\psi(s_f)=0$. 

Now let $A\in \GG^0$ and $b\in B\setminus B_A$. Since $\psi$ is perfect then $\bigcup\limits_{v\in G^0}B_v$ is a basis of $M$ and so there exists a vertex $v\in G^0$ such that $b\in B_v$. Since $b\notin B_A$ then $v\notin A$, and hence $\psi(p_A)\psi(p_v)=0$. Since $b=\psi(p_v)(b)$ we have that $\psi(p_A)(b)=\psi(p_A)\psi(p_v)(b)=0$. \fim

\begin{remark}\label{isomodules} We recall from the proof of Theorem 6.5 in \cite{dd1} that the isomorphism which intertwines the representations $\psi$ and $\varphi$ of the previous proposition is the isomorphism $U:M\rightarrow M(X)$ defined by $U(b_x)=\delta_x$, where $b_x$ is an element of the basis of $M$ and $\delta_x$ is the characteristic map on the set $\{x\}\subseteq X$. For a perfect, permutative representation $\psi:L_R(\GG)\rightarrow End_R(M)$, let $X$ be the perfect branching system as in Proposition \ref{bs1}, and let $T:X\rightarrow \pfin\cup\pinf$ be the morphism of branching systems as in Proposition \ref{morphism}. This morphism $T$ induces a $R$-homomorphism $\Psi:M(X)\rightarrow \mathbb{P}$ defined by $\Psi(\delta_x)=d_{T(x)}$, where $\mathbb{P}$ is the free $R$-module (as in Definition \ref{freemodule}) generated by $\{d_p:p\in \pfin\cup\pinf\}$. So we get a homomorphism $\Phi=\Psi\circ U:M\rightarrow \mathbb{P}$, which takes $b_x \in M$ to $d_{T(x)}$ in $\mathbb{P}$.
\end{remark}

We end the paper characterizing perfect, permutative representations.

\begin{teorema} Let $R$ be a field and $\psi:L_R(\GG)\rightarrow End_R(M)$ be a perfect permutative representation. Then $\psi$ is irreducible if, and only if, $M$ is isomorphic to $\mathbb{P}_{[p]}$ via the homomorphism $\Psi$ (as in the previous remark) for some $p\in \pfin\cup\pinf$.
\end{teorema}


\demo Let $X$ be the branching system of Proposition \ref{bs1} and $\varphi:L_R(\GG)\rightarrow End_R(M(X))$ be the representation induced by this branching system, which is equivalent to $\psi$, by Proposition \ref{eqrep2}. Since $\psi$ and $\varphi$ are equivalent then $\psi$ is irreducible if and only if $\varphi$ is irreducible.

Suppose that $\psi$ is irreducible. Then $\varphi$ is irreducible and it follows, from Theorem~\ref{irreducibletheorem}, that $X$ is isomorphic to $[p]$, for some $p\in \pfin\cup\pinf$, via the morphism $T$ described in the proof of Proposition~ \ref{morphism}. Therefore $\Psi:M(X)\rightarrow \mathbb{P}_{[p]}$ defined on each element $\delta_x$ of the basis of $M(X)$ by $\Psi(\delta_x)=d_{T(x)}$ is an isomorphism and hence $\Phi=\Psi\circ U$ is an isomorphism from $M$ to $\mathbb{P}_{[p]}$

For the converse, suppose that $\Phi:M\rightarrow \mathbb{P}_{[p]}$ is an isomorphism. Then $\Psi:M(X)\rightarrow \mathbb{P}_{[p]}$ is an isomorphism, and hence $T:X\rightarrow [p]$ is a bijection. By the proof of Proposition~\ref{morphism}, $T$ is an isomorphism from the branchig system $X$ to the branching system  $[p]$. By Theorem~\ref{irreducibletheorem} $\varphi$ is irreducible and hence $\psi$ is irreducible.
\fim

\begin{remark}
In the above theorem, the assumption that $R$ is a field is only necessary to show the sufficient condition for irreducibility of $\psi$.
\end{remark}

\begin{exemplo} In the previous theorem, it is important that $M$ is isomorphic to $\mathbb{P}_{[p]}$ via the isomorphism $\Psi$. It is not enough to have $M$ isomorphic to $\mathbb{P}_{[p]}$ via some isomorphism different from $\Psi$. For example, let $\GG$ be a directed graph with infinite edges $\{e,e_1,e_2,e_3,...\}$ and infinite vertices $\{u,w,v_1,v_2,v_3,...\}$ as follows:

\vspace{0.5cm}
\centerline{
\setlength{\unitlength}{1cm}
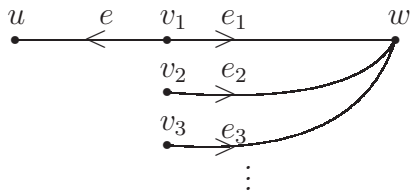
\begin{picture}(0,0)
\put(0,0){\circle*{0.1}}
\put(-0.1,0.2){$v_1$}
\put(-2,0){\circle*{0.1}}
\put(-2.1,0.2){$u$}
\put(-2,0){\line(1,0){2}}
\put(0,0){\line(1,0){3}}
\put(3,0){\circle*{0.1}}
\put(2.9,0.2){$w$}
\put(-1.1,-0.1){$<$}
\put(-0.9,0.2){$e$}
\put(0.6,-0.1){$>$}
\put(0.7,0.2){$e_1$}
\qbezier(0,-0.7)(2.5,-0.9)(3,0)
\qbezier(0,-1.4)(2.5,-1.6)(3,0)
\put(0,-0.7){\circle*{0.1}}
\put(-0.1,-0.5){$v_2$}
\put(0,-1.4){\circle*{0.1}}
\put(-0.1,-1.2){$v_3$}
\put(1,-2){$\vdots$}
\put(0.6,-0.83){$>$}
\put(0.7,-0.5){$e_2$}
\put(0.6,-1.53){$>$}
\put(0.7,-1.3){$e_3$}
\end{picture}}
\vspace{2 cm}

Let $\pi:L_R(\GG)\rightarrow End_R(\mathbb{P})$ be the representation obtained in Proposition \ref{rep1}. This representation is perfect and permutative, following Example \ref{perfectpermutativerep}. Moreover, $\mathbb{P}$ is isomorphic to $\mathbb{P}_{[(w,w)]}$, since both $\mathbb{P}$ and $\mathbb{P}_{[(w,w)]}$ are isomorphic to $\bigoplus\limits_{\N}R$. However $\pi$ is not irreducible since, for example, $\mathbb{P}_{[(e,u)]}$ is $\pi$-invariant.

\end{exemplo}

\vspace{1.5pc}

Daniel Gon\c{c}alves, Departamento de Matem\'{a}tica, Universidade Federal de Santa Catarina, Florian\'{o}polis, 88040-900, Brazil.

Email: daemig@gmail.com

\vspace{0.5pc}
Danilo Royer, Departamento de Matem\'{a}tica, Universidade Federal de Santa Catarina, Florian\'{o}polis, 88040-900, Brazil.

Email: daniloroyer@gmail.com
\vspace{0.5pc}

\end{document}